\documentclass[3p]{elsarticle}
\usepackage{mathtools}
\usepackage{siunitx}
\usepackage{hyperref}
\usepackage[capitalize,nameinlink]{cleveref}
\usepackage[american]{babel}
\usepackage[babel,autostyle=true]{csquotes}
\usepackage{xcolor}
    \definecolor{MatlabBlue}       {rgb}{0.0000,0.4470,0.7410} 
    \definecolor{MatlabRed}        {rgb}{0.8500,0.3250,0.0980}
    \definecolor{MatlabYellow}     {rgb}{0.9290,0.6940,0.1250}
    \definecolor{MatlabPurple}     {rgb}{0.4940,0.1840,0.5560}
    \definecolor{MatlabGreen}      {rgb}{0.4660,0.6740,0.1880}
    \definecolor{MatlabLightblue}  {rgb}{0.3010,0.7450,0.9330}
    \definecolor{MatlabDarkred}    {rgb}{0.6350,0.0780,0.1840}
    \definecolor{MatlabLightred}   {rgb}{1.0000,0.2700,0.2270}
    \definecolor{MatlabCornflower} {rgb}{0.3960,0.5090,0.9920}
    \definecolor{MatlabLightyellow}{rgb}{1.0000,0.8390,0.0390}
    \definecolor{MatlabGold}       {rgb}{0.0000,0.6390,0.6390}
    \definecolor{MatlabBrown}      {rgb}{0.7960,0.5170,0.3640}
\usepackage{subcaption}
\usepackage{placeins}
\usepackage{pgfplots}
    \usetikzlibrary{calc,decorations.markings,decorations.pathmorphing}
    \tikzset{
        house/.pic={
            \draw[fill=#1] (-0.5,0) -- (-0.5,-1) -- (0.5,-1) -- (0.5,0)  -- (0.6,-0.1) -- (0,0.5) -- (-0.6,-0.1);
        },
        tree/.pic={
            \foreach \w/\c in {0.3/30,0.2/50,0.1/70}{
                \fill[
                    brown!\c!black,
                    decoration={random steps, segment length=2, amplitude=0.2},
                    decorate,
                    transform shape,
                ] (0,0) ++(-\w/2,0) rectangle +(\w,-3);
            }
            \foreach \n/\c in {1.4/40,1.2/50,1/60,0.8/70,0.6/80,0.4/90} {
               \fill[
                    green!\c!black,
                    decoration={random steps, segment length=0.4, amplitude=2},
                    decorate,
                    transform shape,
                ] ellipse (\n/1.5 and \n);
            }%
        },
        network edge/.style={
            draw,
            decoration={
                markings,
                mark=at position 0.35 with {\arrow{latex}},
                mark=at position 0.85 with {\arrow{latex}}
            },
            postaction={decorate},
        },
        network node/.style={
            draw,
            circle,
            minimum width=14pt,
            inner sep=0pt,
            outer sep=0pt,
        },
    }
    \pgfplotscreateplotcyclelist{MatlabColors}{
        {MatlabBlue},
        {MatlabRed},
        {MatlabYellow},
        {MatlabPurple},
        {MatlabGreen},
        {MatlabLightblue},
        {MatlabDarkred},
        {MatlabLightred},
        {MatlabCornflower},
        {MatlabLightyellow},
        {MatlabGold},
        {MatlabBrown}
    }
    \pgfplotsset{
        compat=1.18,
        every axis/.style={
            legend style={font=\footnotesize},
            legend cell align={left},
            cycle list name=MatlabColors,
        },
        every axis plot/.append style={thick},
        /pgf/number format/1000 sep = \thinspace,
        colormap={turbo}{
            [1pt] rgb(0pt)=(0.18995,0.07176,0.23217); rgb(1pt)=(0.19483,0.08339,0.26149); rgb(2pt)=(0.19956,0.09498,0.29024); rgb(3pt)=(0.20415,0.10652,0.31844); rgb(4pt)=(0.2086,0.11802,0.34607); rgb(5pt)=(0.21291,0.12947,0.37314); rgb(6pt)=(0.21708,0.14087,0.39964); rgb(7pt)=(0.22111,0.15223,0.42558); rgb(8pt)=(0.225,0.16354,0.45096); rgb(9pt)=(0.22875,0.17481,0.47578); rgb(10pt)=(0.23236,0.18603,0.50004); rgb(11pt)=(0.23582,0.1972,0.52373); rgb(12pt)=(0.23915,0.20833,0.54686); rgb(13pt)=(0.24234,0.21941,0.56942); rgb(14pt)=(0.24539,0.23044,0.59142); rgb(15pt)=(0.2483,0.24143,0.61286); rgb(16pt)=(0.25107,0.25237,0.63374); rgb(17pt)=(0.25369,0.26327,0.65406); rgb(18pt)=(0.25618,0.27412,0.67381); rgb(19pt)=(0.25853,0.28492,0.693); rgb(20pt)=(0.26074,0.29568,0.71162); rgb(21pt)=(0.2628,0.30639,0.72968); rgb(22pt)=(0.26473,0.31706,0.74718); rgb(23pt)=(0.26652,0.32768,0.76412); rgb(24pt)=(0.26816,0.33825,0.7805); rgb(25pt)=(0.26967,0.34878,0.79631); rgb(26pt)=(0.27103,0.35926,0.81156); rgb(27pt)=(0.27226,0.3697,0.82624); rgb(28pt)=(0.27334,0.38008,0.84037); rgb(29pt)=(0.27429,0.39043,0.85393); rgb(30pt)=(0.27509,0.40072,0.86692); rgb(31pt)=(0.27576,0.41097,0.87936); rgb(32pt)=(0.27628,0.42118,0.89123); rgb(33pt)=(0.27667,0.43134,0.90254); rgb(34pt)=(0.27691,0.44145,0.91328); rgb(35pt)=(0.27701,0.45152,0.92347); rgb(36pt)=(0.27698,0.46153,0.93309); rgb(37pt)=(0.2768,0.47151,0.94214); rgb(38pt)=(0.27648,0.48144,0.95064); rgb(39pt)=(0.27603,0.49132,0.95857); rgb(40pt)=(0.27543,0.50115,0.96594); rgb(41pt)=(0.27469,0.51094,0.97275); rgb(42pt)=(0.27381,0.52069,0.97899); rgb(43pt)=(0.27273,0.5304,0.98461); rgb(44pt)=(0.27106,0.54015,0.9893); rgb(45pt)=(0.26878,0.54995,0.99303); rgb(46pt)=(0.26592,0.55979,0.99583); rgb(47pt)=(0.26252,0.56967,0.99773); rgb(48pt)=(0.25862,0.57958,0.99876); rgb(49pt)=(0.25425,0.5895,0.99896); rgb(50pt)=(0.24946,0.59943,0.99835); rgb(51pt)=(0.24427,0.60937,0.99697); rgb(52pt)=(0.23874,0.61931,0.99485); rgb(53pt)=(0.23288,0.62923,0.99202); rgb(54pt)=(0.22676,0.63913,0.98851); rgb(55pt)=(0.22039,0.64901,0.98436); rgb(56pt)=(0.21382,0.65886,0.97959); rgb(57pt)=(0.20708,0.66866,0.97423); rgb(58pt)=(0.20021,0.67842,0.96833); rgb(59pt)=(0.19326,0.68812,0.9619); rgb(60pt)=(0.18625,0.69775,0.95498); rgb(61pt)=(0.17923,0.70732,0.94761); rgb(62pt)=(0.17223,0.7168,0.93981); rgb(63pt)=(0.16529,0.7262,0.93161); rgb(64pt)=(0.15844,0.73551,0.92305); rgb(65pt)=(0.15173,0.74472,0.91416); rgb(66pt)=(0.14519,0.75381,0.90496); rgb(67pt)=(0.13886,0.76279,0.8955); rgb(68pt)=(0.13278,0.77165,0.8858); rgb(69pt)=(0.12698,0.78037,0.8759); rgb(70pt)=(0.12151,0.78896,0.86581); rgb(71pt)=(0.11639,0.7974,0.85559); rgb(72pt)=(0.11167,0.80569,0.84525); rgb(73pt)=(0.10738,0.81381,0.83484); rgb(74pt)=(0.10357,0.82177,0.82437); rgb(75pt)=(0.10026,0.82955,0.81389); rgb(76pt)=(0.0975,0.83714,0.80342); rgb(77pt)=(0.09532,0.84455,0.79299); rgb(78pt)=(0.09377,0.85175,0.78264); rgb(79pt)=(0.09287,0.85875,0.7724); rgb(80pt)=(0.09267,0.86554,0.7623); rgb(81pt)=(0.0932,0.87211,0.75237); rgb(82pt)=(0.09451,0.87844,0.74265); rgb(83pt)=(0.09662,0.88454,0.73316); rgb(84pt)=(0.09958,0.8904,0.72393); rgb(85pt)=(0.10342,0.896,0.715); rgb(86pt)=(0.10815,0.90142,0.70599); rgb(87pt)=(0.11374,0.90673,0.69651); rgb(88pt)=(0.12014,0.91193,0.6866); rgb(89pt)=(0.12733,0.91701,0.67627); rgb(90pt)=(0.13526,0.92197,0.66556); rgb(91pt)=(0.14391,0.9268,0.65448); rgb(92pt)=(0.15323,0.93151,0.64308); rgb(93pt)=(0.16319,0.93609,0.63137); rgb(94pt)=(0.17377,0.94053,0.61938); rgb(95pt)=(0.18491,0.94484,0.60713); rgb(96pt)=(0.19659,0.94901,0.59466); rgb(97pt)=(0.20877,0.95304,0.58199); rgb(98pt)=(0.22142,0.95692,0.56914); rgb(99pt)=(0.23449,0.96065,0.55614); rgb(100pt)=(0.24797,0.96423,0.54303); rgb(101pt)=(0.2618,0.96765,0.52981); rgb(102pt)=(0.27597,0.97092,0.51653); rgb(103pt)=(0.29042,0.97403,0.50321); rgb(104pt)=(0.30513,0.97697,0.48987); rgb(105pt)=(0.32006,0.97974,0.47654); rgb(106pt)=(0.33517,0.98234,0.46325); rgb(107pt)=(0.35043,0.98477,0.45002); rgb(108pt)=(0.36581,0.98702,0.43688); rgb(109pt)=(0.38127,0.98909,0.42386); rgb(110pt)=(0.39678,0.99098,0.41098); rgb(111pt)=(0.41229,0.99268,0.39826); rgb(112pt)=(0.42778,0.99419,0.38575); rgb(113pt)=(0.44321,0.99551,0.37345); rgb(114pt)=(0.45854,0.99663,0.3614); rgb(115pt)=(0.47375,0.99755,0.34963); rgb(116pt)=(0.48879,0.99828,0.33816); rgb(117pt)=(0.50362,0.99879,0.32701); rgb(118pt)=(0.51822,0.9991,0.31622); rgb(119pt)=(0.53255,0.99919,0.30581); rgb(120pt)=(0.54658,0.99907,0.29581); rgb(121pt)=(0.56026,0.99873,0.28623); rgb(122pt)=(0.57357,0.99817,0.27712); rgb(123pt)=(0.58646,0.99739,0.26849); rgb(124pt)=(0.59891,0.99638,0.26038); rgb(125pt)=(0.61088,0.99514,0.2528); rgb(126pt)=(0.62233,0.99366,0.24579); rgb(127pt)=(0.63323,0.99195,0.23937); rgb(128pt)=(0.64362,0.98999,0.23356); rgb(129pt)=(0.65394,0.98775,0.22835); rgb(130pt)=(0.66428,0.98524,0.2237); rgb(131pt)=(0.67462,0.98246,0.2196); rgb(132pt)=(0.68494,0.97941,0.21602); rgb(133pt)=(0.69525,0.9761,0.21294); rgb(134pt)=(0.70553,0.97255,0.21032); rgb(135pt)=(0.71577,0.96875,0.20815); rgb(136pt)=(0.72596,0.9647,0.2064); rgb(137pt)=(0.7361,0.96043,0.20504); rgb(138pt)=(0.74617,0.95593,0.20406); rgb(139pt)=(0.75617,0.95121,0.20343); rgb(140pt)=(0.76608,0.94627,0.20311); rgb(141pt)=(0.77591,0.94113,0.2031); rgb(142pt)=(0.78563,0.93579,0.20336); rgb(143pt)=(0.79524,0.93025,0.20386); rgb(144pt)=(0.80473,0.92452,0.20459); rgb(145pt)=(0.8141,0.91861,0.20552); rgb(146pt)=(0.82333,0.91253,0.20663); rgb(147pt)=(0.83241,0.90627,0.20788); rgb(148pt)=(0.84133,0.89986,0.20926); rgb(149pt)=(0.8501,0.89328,0.21074); rgb(150pt)=(0.85868,0.88655,0.2123); rgb(151pt)=(0.86709,0.87968,0.21391); rgb(152pt)=(0.8753,0.87267,0.21555); rgb(153pt)=(0.88331,0.86553,0.21719); rgb(154pt)=(0.89112,0.85826,0.2188); rgb(155pt)=(0.8987,0.85087,0.22038); rgb(156pt)=(0.90605,0.84337,0.22188); rgb(157pt)=(0.91317,0.83576,0.22328); rgb(158pt)=(0.92004,0.82806,0.22456); rgb(159pt)=(0.92666,0.82025,0.2257); rgb(160pt)=(0.93301,0.81236,0.22667); rgb(161pt)=(0.93909,0.80439,0.22744); rgb(162pt)=(0.94489,0.79634,0.228); rgb(163pt)=(0.95039,0.78823,0.22831); rgb(164pt)=(0.9556,0.78005,0.22836); rgb(165pt)=(0.96049,0.77181,0.22811); rgb(166pt)=(0.96507,0.76352,0.22754); rgb(167pt)=(0.96931,0.75519,0.22663); rgb(168pt)=(0.97323,0.74682,0.22536); rgb(169pt)=(0.97679,0.73842,0.22369); rgb(170pt)=(0.98,0.73,0.22161); rgb(171pt)=(0.98289,0.7214,0.21918); rgb(172pt)=(0.98549,0.7125,0.2165); rgb(173pt)=(0.98781,0.7033,0.21358); rgb(174pt)=(0.98986,0.69382,0.21043); rgb(175pt)=(0.99163,0.68408,0.20706); rgb(176pt)=(0.99314,0.67408,0.20348); rgb(177pt)=(0.99438,0.66386,0.19971); rgb(178pt)=(0.99535,0.65341,0.19577); rgb(179pt)=(0.99607,0.64277,0.19165); rgb(180pt)=(0.99654,0.63193,0.18738); rgb(181pt)=(0.99675,0.62093,0.18297); rgb(182pt)=(0.99672,0.60977,0.17842); rgb(183pt)=(0.99644,0.59846,0.17376); rgb(184pt)=(0.99593,0.58703,0.16899); rgb(185pt)=(0.99517,0.57549,0.16412); rgb(186pt)=(0.99419,0.56386,0.15918); rgb(187pt)=(0.99297,0.55214,0.15417); rgb(188pt)=(0.99153,0.54036,0.1491); rgb(189pt)=(0.98987,0.52854,0.14398); rgb(190pt)=(0.98799,0.51667,0.13883); rgb(191pt)=(0.9859,0.50479,0.13367); rgb(192pt)=(0.9836,0.49291,0.12849); rgb(193pt)=(0.98108,0.48104,0.12332); rgb(194pt)=(0.97837,0.4692,0.11817); rgb(195pt)=(0.97545,0.4574,0.11305); rgb(196pt)=(0.97234,0.44565,0.10797); rgb(197pt)=(0.96904,0.43399,0.10294); rgb(198pt)=(0.96555,0.42241,0.09798); rgb(199pt)=(0.96187,0.41093,0.0931); rgb(200pt)=(0.95801,0.39958,0.08831); rgb(201pt)=(0.95398,0.38836,0.08362); rgb(202pt)=(0.94977,0.37729,0.07905); rgb(203pt)=(0.94538,0.36638,0.07461); rgb(204pt)=(0.94084,0.35566,0.07031); rgb(205pt)=(0.93612,0.34513,0.06616); rgb(206pt)=(0.93125,0.33482,0.06218); rgb(207pt)=(0.92623,0.32473,0.05837); rgb(208pt)=(0.92105,0.31489,0.05475); rgb(209pt)=(0.91572,0.3053,0.05134); rgb(210pt)=(0.91024,0.29599,0.04814); rgb(211pt)=(0.90463,0.28696,0.04516); rgb(212pt)=(0.89888,0.27824,0.04243); rgb(213pt)=(0.89298,0.26981,0.03993); rgb(214pt)=(0.88691,0.26152,0.03753); rgb(215pt)=(0.88066,0.25334,0.03521); rgb(216pt)=(0.87422,0.24526,0.03297); rgb(217pt)=(0.8676,0.2373,0.03082); rgb(218pt)=(0.86079,0.22945,0.02875); rgb(219pt)=(0.8538,0.2217,0.02677); rgb(220pt)=(0.84662,0.21407,0.02487); rgb(221pt)=(0.83926,0.20654,0.02305); rgb(222pt)=(0.83172,0.19912,0.02131); rgb(223pt)=(0.82399,0.19182,0.01966); rgb(224pt)=(0.81608,0.18462,0.01809); rgb(225pt)=(0.80799,0.17753,0.0166); rgb(226pt)=(0.79971,0.17055,0.0152); rgb(227pt)=(0.79125,0.16368,0.01387); rgb(228pt)=(0.7826,0.15693,0.01264); rgb(229pt)=(0.77377,0.15028,0.01148); rgb(230pt)=(0.76476,0.14374,0.01041); rgb(231pt)=(0.75556,0.13731,0.00942); rgb(232pt)=(0.74617,0.13098,0.00851); rgb(233pt)=(0.73661,0.12477,0.00769); rgb(234pt)=(0.72686,0.11867,0.00695); rgb(235pt)=(0.71692,0.11268,0.00629); rgb(236pt)=(0.7068,0.1068,0.00571); rgb(237pt)=(0.6965,0.10102,0.00522); rgb(238pt)=(0.68602,0.09536,0.00481); rgb(239pt)=(0.67535,0.0898,0.00449); rgb(240pt)=(0.66449,0.08436,0.00424); rgb(241pt)=(0.65345,0.07902,0.00408); rgb(242pt)=(0.64223,0.0738,0.00401); rgb(243pt)=(0.63082,0.06868,0.00401); rgb(244pt)=(0.61923,0.06367,0.0041); rgb(245pt)=(0.60746,0.05878,0.00427); rgb(246pt)=(0.5955,0.05399,0.00453); rgb(247pt)=(0.58336,0.04931,0.00486); rgb(248pt)=(0.57103,0.04474,0.00529); rgb(249pt)=(0.55852,0.04028,0.00579); rgb(250pt)=(0.54583,0.03593,0.00638); rgb(251pt)=(0.53295,0.03169,0.00705); rgb(252pt)=(0.51989,0.02756,0.0078); rgb(253pt)=(0.50664,0.02354,0.00863); rgb(254pt)=(0.49321,0.01963,0.00955); rgb(255pt)=(0.4796,0.01583,0.01055)
      },
   }

\usepackage{xspace}
\newcommand{\ie}[0]{i.e.\@\xspace}
\newcommand{\eg}[0]{e.g.\@\xspace}
\newcommand{\nbhy}[0]{\babelhyphen{nobreak}}

\newcommand*{\dd}{\mathrm{d}}            
\newcommand*{\tp}{\mathrm{T}}            
\newcommand*{\loc}{\mathrm{loc}}         
\newcommand*{\IN}[0]{\mathbb{N}}         
\newcommand*{\IR}[0]{\mathbb{R}}         
\newcommand*{\pder}[1]{\tfrac{\partial}{\partial #1}} 
\renewcommand{\vec}[1]{\boldsymbol{#1}}  
\makeatletter
\newcommand*{\Cases}[4][]{
    \bgroup%
    \newcommand{\myEndDelim}{#1}
    \newcommand{\myInnerDelim}[0]{,}
    \newcommand{\myOuterDelim}[0]{,}
    \newcommand{\inCase}[0]{\text{if }}
    \newcommand{\otherwise}[0]{\text{else}}
    \begin{cases}
        #2 \myInnerDelim & \inCase#3\myOuterDelim \\
        #4 \myInnerDelim & \checkNextRightArg
    }
    \newcommand{\checkNextRightArg}[0]{\@ifnextchar\bgroup{\gobbleNextRightArg}{\otherwise \myEndDelim \end{cases}\egroup}}
    \newcommand{\checkNextLeftArg}[0]{\@ifnextchar\bgroup{\myOuterDelim \\ \gobbleNextLeftArg}{\myEndDelim \end{cases}\egroup}}
    \newcommand{\gobbleNextRightArg}[1]{\inCase #1 \checkNextLeftArg}
    \newcommand{\gobbleNextLeftArg}[1]{#1 \myInnerDelim & \checkNextRightArg}
\makeatother

\newcommand*{\ur}[1][]{#1q^\mathrm{r}}        
\newcommand*{\ul}[1][]{#1q^\mathrm{\ell}}     
\newcommand{\influx}[1][]{#1\phi^{\mathrm{in}}} 
\newcommand{\outflux}[1][]{#1\phi^{\mathrm{out}}} 
\newcommand{\bmax}[0]{b^{\mathrm{max}}}  
\newcommand{\qmax}[0]{q^{\mathrm{max}}}  
\newcommand{\trt}[0]{\daleth}            
\newcommand{\NT}[0]{\mathrm{NT}}         

\makeatletter
\let\oldfigure\figure
\def\figure{\@ifnextchar[\figure@i \figure@ii}
\def\figure@i[#1]{\oldfigure[#1]\centering} 
\def\figure@ii{\oldfigure\centering}        

\let\oldsubfigure\subfigure
\def\subfigure{\@ifnextchar[\subfigure@i \subfigure@ii}
\def\subfigure@i[#1]#2{\hfil\oldsubfigure[#1]{#2}\centering} 
\def\subfigure@ii#1{\hfil\oldsubfigure{#1}\centering}        
\let\oldendsubfigure\endsubfigure
\def\endsubfigure{\oldendsubfigure\hfil\mbox{}}

\let\oldtable\table
\def\table{\@ifnextchar[\table@i \table@ii}
\def\table@i[#1]{\oldtable[#1]\centering} 
\def\table@ii{\oldtable\centering}        

\usepackage{makecell}
\usepackage{autonum}
\newcommand{\emphTwo}[1]{\textit{#1}}
\newcommand*{\loadandscalepgfplot}[2]{%
   \bgroup%
   \pgfplotsset{
      every axis/.append style={scale=#1},    
      every colorbar/.append style={scale=1}, 
   }
   \input{#2}%
   \egroup
}

\usepackage{graphicx}%
\usepackage{multirow}%
\usepackage{amsmath,amssymb,amsfonts}%
\usepackage{amsthm}%
\usepackage{mathrsfs}%
\usepackage{lmodern}

\usepackage{xcolor}%
\usepackage{booktabs}%
\usepackage{algorithm}%
\usepackage{algorithmicx}%
\usepackage{algpseudocode}%
\usepackage{orcidlink}

\newtheorem{remark}{Remark}%

\newtheorem{definition}{Definition}%

\begin{document}
\title{A general framework for nonlocal traffic flow models on networks}

\author[1]{Alexander Keimer\,\orcidlink{0000-0003-3825-5853}}\ead{alexander.keimer@uni-rostock.de}

\author[2,3]{Lukas Pflug\,\orcidlink{0000-0001-8001-5832}}\ead{lukas.pflug@fau.de}

\author[2,3]{Florian Prohaska\,\orcidlink{0009-0004-8034-6366}}\ead{florian.prohaska@fau.de}

\affiliation[1]{organisation={Institute of Mathematics, University of Rostock}, address={Ulmenstraße~69}, city={Rostock}, postcode={18057}, country={Germany}}

\affiliation[2]{organisation={FAU Competence Center Scientific Computing, Friedrich-Alexander-University Erlangen-Nürnberg (FAU)}, address={Martensstraße~5a}, city={Erlangen}, postcode={91058}, country={Germany}}

\affiliation[3]{organisation={Chair of Applied Mathematics (Continuous Optimisation), Friedrich-Alexander-University Erlangen-Nürnberg (FAU)}, address={Cauerstraße~11}, city={Erlangen}, postcode={91058}, country={Germany}}

\begin{abstract}We present a general computational framework for macroscopic nonlocal traffic flow models on networks with multiple commodities. The model combines scalar conservation laws on edges with nonlocal velocity functions and couples them via buffer-based junction dynamics.

Routing is defined in general as a prescription for distributing drivers across outgoing road segments. As an example, we implement dynamic k-shortest-path routing, where travel times along roads and waiting times at intersections are used to compute shortest paths to the commodities' destinations at each time step, and drivers are distributed accordingly. In another example, we optimize routing over a considered time horizon to minimize the total travel time. This framework naturally creates a feedback loop between traffic evolution and route choice. 

Numerical examples, ranging from small test cases to large grid-like networks, demonstrate the robustness of the approach and allow for a comparison of different routing strategies.
\end{abstract}


\maketitle

\section{Introduction}
Although more and more real-time traffic data is becoming available through, \eg, routing apps, connected cars, and autonomous driving systems, their potential for traffic optimization is far from being explored exhaustively. While single roads might experience traffic flow smoothing \cite{Lee2025} by taking specific measures, optimizing traffic flow on the city network scale is not intensively studied when considering time dependent routing. Proper objectives to optimize could be network performance, averaged travel time, emissions, and more, while control would mainly consist of changing routing at the intersections over time. To address the task of significantly improving traffic performance on the macroscopic level, we propose a scalable macroscopic traffic flow model for general networks and several ideas of routing algorithms.

On the modelling level, we represent each road of a given network by a directed link on which we define a system of nonlocal conservation laws. The intersections are modelled by bounded queues, taking into account out and inflow, routing suggestions and priorization of incoming roads.

The task, we thus have in mind is also sometimes referred to as dynamic traffic assignment (DTA) \cite{DTA_primer}. Most models coming up in this research field~\cite{Keimer2020Routing,Bayen2019Time,Papageorgiou1990Dynamic,Saw2015Literature} specify, just as we do, traffic dynamics on each road, junction models, and routing for origin--destination pairs.

\cite{garavello2006traffic,garavello2016models} provides comprehensive results for modeling traffic flow with conservation laws on networks. Concerning the existence and uniqueness of entropy solutions on networks for local conservation laws and as such local traffic flow models like the LWR (after Lighthill, Witham, and Richards \cite{lwr_1,lwr_2}) as well as the proper definitions, one finds substantial results in \cite{Musch2022,Garavello2009,Bressan2015,Ancona2018,LaurentBrouty2020,Samaranayake2018}.
Concerning related gas/traffic dynamics on networks, one finds results in e.g., 
\cite{Kotsialos2002,Egger2022,Egger2017Numerical,Bongarti2024,Samaranayake2018,Festa2023}; however, only a prototypical implementation of corresponding macroscopic models is given, or as in \cite{FOKKEN2022}, a fast numerical solver is provided  but for gas networks with no routing or nonlocality. In \cite{Leclercq2015} conservation laws on networks are used, but link dynamics are entirely replaced by single origin--destination conservation laws and the demand is realized as a source term. \cite{Thonhofer2018} is more in the spirit of what we plan to model, but once more for local conservation laws, merely on the discretized level and not considering dynamic routing.

There are also many commercial toolboxes for parts of the dynamics traffic assignment problem available, we only name a few: SUMO \cite{behrisch2011sumo}, \href{https://www.anylogic.com/road-traffic/}{Anylogic}, \href{https://www.ptvgroup.com/en/products/ptv-vissim}{PTV-VISSIM}, Cityflow \cite{Tang2019}, \href{https://www.aimsun.com/}{Aimsun}, and MATSim \cite{horni2016introducing}. However, these toolboxes use mostly microscopic models, which are not scalable as the computational complexity increases rapidly with the volume of traffic and they are also difficult to use when it comes to tasks beyond their provided capabilities (optimization via adjoints, etc.).

When it comes to routing choices, instantaneous and past time routing as well as game theoretical decision making on the continuous level \cite{von2007theory,bacsar1998dynamic,Wishart1966}, we refer the reader to \cite{Patriksson2015TheTA} for stating the definition of Wardrop's first and second principle (the traffic analogous to game theoretical equilibria concerning routing) and to \cite{Bayen2019Time,Keimer2020Routing,Gao2010,Bovy2009,BenAkiva2004,Engelson2003,Chai2017,Graf2020,Graf2023Prediction,Graf2025Dynamic} and the references therein together with the disclaimer that such a list cannot be complete due to the huge amount of publications in this field.

In \cref{sec:model}, we introduce the model for a network of streets, including the nonlocal impact, the buffer dynamics and define routing in a general manner, supported by an exemplary routing algorithm based on travel times.
\Cref{sec:notes} provides some notes on the implementation used to compute the examples.
The latter are collected \cref{sec:examples} and include a graspable introductory example focusing on expected routing behavior, a New York-like street network demonstrating a nontrivial example, and finally an academic example to explain the quality of different routing algorithms.
\section{Model} \label{sec:model}
\begin{definition}[Basic model of a street]
    Consider a street of length~$L$, on which we impose the initial initial-boundary value problem~\cite{Keimer2018Nonlocal}
    \begin{alignat}{2}
        \pder{t}q(t,x) + \pder{x}\bigl(q(t,x) v[q](t,x)\bigr) &= 0               , & \quad (t,x) &\in (0,T)\times(0,L)\eqqcolon\Omega_T, \\
        q(0,x)                                                &= q_0(x)          , & \quad    x  &\in (0,L),   \label{eq:initialDensity} \\
        q(t,0) v[q](t,0)                                      &= \ul(t) v[q](t,0), & \quad  t    &\in (0,T),
    \end{alignat}
    with time horizon $T\in\IR_{>0}$, density $q\colon\Omega_T\to\IR_{\geq0}$, velocity $v[q]\colon\Omega_T\to\IR_{\geq0}$, initial density $q_0\colon(0,L)\to\IR_{\geq0}$ at time $t=0$, and boundary datum $\ul\colon(0,T)\to\IR_{\geq0}$ at the beginning of the street.
\end{definition}
The velocity function~$v[q]$ is supplemented by the nonlocal impact~$W$ in the following sense~\cite{Keimer2019approximation,Keimer2018Nonlocal}.
\begin{definition}[Nonlocal impact] \label{def:nonlocalImpact}
    We define
    \begin{align}
        \label{eq:velocity}
        v[q]\colon \Omega_T\ni(t,x) &\mapsto V\bigl(W[q, \gamma, \eta](t,x)\bigr), \\
        \label{eq:nonlocalImpact}
        W[q, \gamma, \eta]\colon \Omega_T\ni(t,x)
            &\mapsto
            \tfrac{1}{\eta} \biggl(
                       \int_x^L      \gamma\bigl(\tfrac{y-x}{\eta}\bigr) q(t, y) \,\dd y + 
                \ur(t) \int_L^\infty \gamma\bigl(\tfrac{y-x}{\eta}\bigr)         \,\dd y
            \biggr),
    \end{align}
    where $V\in W^{1,\infty}_\loc(\IR_{\geq0};\IR_{\geq0})$ is non-increasing and fulfills $V(\qmax)=0$ for the maximum density~$\qmax$ on the street, $\eta\in\IR_{>0}$ the look-ahead distance, $\gamma\colon\IR_{\geq0}\to\IR_{\geq0}$ a non-increasing kernel with unit mass and $\ur$~is the boundary datum at the right end of the street.
\end{definition}
As long as the nonlocal term \enquote{sees} the street, we convolve the kernel with the density present, while the \enquote{rest} of the kernel is convolved with the right boundary datum~$\ur$. Consequently, the nonlocal term evaluated at the end of the street, \ie, at $x=L$, yields the right boundary datum, \ie, $W[q,\gamma,\eta](t, L) = \ur(t)$, $t\in(0,T)$.

In preparation of developing proper dynamics on the network level, we introduce buffers~\cite{LaurentBrouty2020} at the beginning of each street which model the interactions between roads connected by a junction.
\begin{definition}[Buffer dynamics] \label{def:bufferDynamics}
    The time-dependent buffer load $b\colon(0,T)\to\IR_{\geq0}$ is determined by the initial value problem
    \begin{alignat}{2}
        \dot b(t) &= \influx(t) - \outflux(t), & \quad t &\in (0,T), \\
        b(0)      &= b_0, \label{eq:initialBufferLoad}
    \end{alignat}
    where $b_0\in\IR_{\geq0}$ denotes the initial buffer load, and $\influx\colon(0,T)\to\IR_{\geq0}$ and $\outflux\colon(0,T)\to\IR_{\geq0}$ the buffer's influx and outflux, respectively.
\end{definition}
There are several ways to place the named buffers in the network, each of them having different advantages and disadvantages. Throughout this work, we pose one and only one buffer at the beginning of each street (\ie, we do not keep track of all the roads entering the buffer separately), controlling the flux onto the associated road and collecting the mass from all preceding streets.

Further, we introduce commodities in order to be able to model different groups of drivers heading to different destinations. Those commodities could also represent different types of traffic participants. In addition, when it comes to routing, we could incorporate preference of one commodity over another. To distinguish between them, we use the following notation.

\begin{remark}[Notation and conventions]
    \begin{enumerate}
        \item We extend the density, buffer load, in- and outflux to be vector valued functions~$\vec q$, $\vec b$, $\influx[\vec]$ and~$\outflux[\vec]$, respectively, where each entry corresponds to one commodity. The scalar quantities~$q$, $b$, $\influx$, and~$\outflux$ correspond to the sum over the individual commodities. If not stated differently, products of vectors are to be understood as element-wise.
        \item Since all of the above quantities may already have subindices indicating the corresponding street, we refer to the entry for the $c$\nbhy th commodity by $[\ast]_c=\vec e_c^\tp\cdot\ast$, \ie, the dot product of the $c$\nbhy th unit vector with the quantity vector.
        \item Furthermore, denote the set of successors of a street~$s$ by~$\zeta_s^+$, where a successor is a street starting at the end of~$s$. Analogously, define~$\zeta_s^-$ to be the set of predecessors of~$s$.
    \end{enumerate}
\end{remark}

We still owe formulae for a number of quantities. To this end, we first need to introduce the vector-valued routing function~\cite{Keimer2020Routing} in a very general manner, which is essential for extending traffic flow models from single roads or intersections to city-scale networks.
\begin{definition}[Routing] \label{def:routing}
    Having the total mass which leaves street~$s$ at time $t\in(0,T)$ and belongs to the $c$\nbhy th commodity, the $c$\nbhy th component of $\vec p_{s,r}(t)\geq\boldsymbol0$ gives the share of mass heading to street~$r$. Consequently, for each street~$s$ which is not a sink---\ie, $\lvert\zeta_s^+\lvert>0$---and from which the destination can be reached, $\sum_r \vec p_{s,r}(t) = \vec 1$ must hold for each time~$t$. In addition, $\vec p_{s,r}(t)=\vec 0$ if $s\notin\zeta_r^-$. A collection of such functions for all combinations of streets~$s$ and~$r$ is said to be a routing.
\end{definition}
Choosing a meaningful and realistic routing is a delicate topic as it describes human behavior. Ignoring fixed splitting ratios, instantaneous $k$~shortest path routing (\cref{def:kShortestPath}) is among the simplest routing rules. In the most complex case, routing could be model-predictive and could therefore be subject to an optimal control type of problem.

Recalling the buffer dynamics in \cref{def:bufferDynamics} we specify the in- and outflux of the buffer, which can in general be functionals of the entire traffic state of the network. However, for the ease of notation, we choose rather straightforward dependencies.
\begin{definition}[Buffer influx] \label{def:bufferInflux}
    For a street~$s$ the influx to the buffer at time~$t\in(0,T)$ is given by
    \begin{equation}
        \influx[\vec]_s(t) = \sum_{r\in\zeta_s^-} \vec p_{r,s}(t) \vec q_r(t, L_r) v_r[q_r](t, L_r),
    \end{equation}
    where the subscripts~$s$ and~$r$ mark properties belonging to the respective streets and $v_r[q_r]$ denotes the nonlocal velocity from \cref{eq:velocity}.
\end{definition}

\begin{definition}[Right boundary datum]
    The right boundary datum is given by the maximum relative buffer load that is seen by any driver who is currently close to the end of their street, rescaled to the maximum density on the street, \ie,
\begin{equation}
    \ur_s(t) \coloneqq \qmax_s \max_{\{ r\in\zeta_s^+ \,\vert\, \exists c\colon [\vec q_s(t, L_s)]_c > 0 \land [\vec p_{s, r}(t)]_c > 0 \}}
   \tfrac{b_r(t)}{\bmax_r}
\end{equation}
for $t\in(0,T)$ and maximum buffer load $\bmax\in\IR_{\geq0}$ ($\max\emptyset\coloneqq0$).
\end{definition}

Given a left boundary datum~$\ul[\vec]$, we can define the outflux of a buffer.
\begin{definition}[Buffer outflux]
    The flux from a buffer to its associated street~$s$ is given by
    \begin{equation}
        \outflux[\vec]_s(t) \coloneqq \ul[\vec]_s(t) v_s[q_s](t, 0)
    \end{equation}
    with $\ul[\vec]_s$ such that~$\ul_s$ is a monotonously increasing with increasing~$\vec b_s(t)$. Further, whenever a component of~$\vec b(t)$ is~$0$, so is the same component of $\ul[\vec]_s(t)$.
\end{definition}
From a modeling perspective, it appears sensible to choose~$\ul_s$ to be a sufficiently smooth lower approximation of a function that is~$0$ when $b_s(t)$ is, and $\qmax_s$ otherwise.

In contrast to the influx of a buffer (\cref{def:bufferInflux}), there is no natural choice for the left boundary datum governing the flux from a buffer to its associated street. Therefore, we only pose the above mild conditions on the function~$\ul[\vec]_s$. We demand $[\vec a]_c=0$ to imply $[\ul[\vec]_s(t)]_c=0$ to ensure that no mass can flow onto the street if the buffer share of the corresponding commodity is empty. There are choices of~$\ul[\vec]_s$ ensuring buffer evacuation in finite time in case of $\influx_s\equiv0$, including
\begin{equation} \label{eq:influxDensity}
    \ul[\vec]_s \colon t \longmapsto \Cases{
        \Bigl(\frac{\vec b_s(t)}{\bmax_s}\Bigr)^{\frac{1}{m_s}} \qmax_s \frac{\vec b_s(t)}{b_s(t)}
    }{
        b_s(t)>0
    }{
        \boldsymbol0
    }
\end{equation}
for $m_s\geq2$. The factor $(\frac{\vec b_s(t)}{\bmax_s(t)})^{\frac{1}{m_s}}$ can be seen element-wise as a softened sign function ensuring continuity. The aforementioned instance \cref{eq:influxDensity} of~$\ul[\vec]_s$ is used throughout the numerical examples in \cref{sec:examples}. Not only it incorporates all the desirable properties, but also allows the network to evacuate in finite time in case of no flux entering the network.

So far, we have not specified potential instantiations of routing, cf.\@ \cref{def:routing}. One possibility is instantaneous $k$~shortest path routing for $k\in\IN_{\geq1}$: Given the travel time~$\trt_s(t)$ of all streets~$s$ at time~$t\in(0,T)$, \ie, the time required to travel from the start to the end of the respective street and potentially the waiting time in the buffer---which can be evaluated (almost instantaneously) using the computed density---we can compute for each node and each destination of a commodity the $k$~shortest paths. Based on the time required to pass these paths, the routing~$\vec p$ is assigned.
\begin{definition}[$k$~shortest path routing] \label{def:kShortestPath}
    Fix a street~$s$ and a commodity~$c$. Assume a $k$~shortest path algorithm found $\ell\leq k$ paths $\wp_r$ with travel times~$\trt_{\wp_r}$ from the end of street~$s$ to the commodity's destination, $r=1,\ldots,\ell$. If~$s$ already leads to the destination, set the $c$\nbhy th component of $\vec p_{s,\hat s}$ to~0 for each street~$\hat s$. Otherwise use a nonnegative non-increasing weighting function~$g$ to compute
    \begin{align} \label{eq:pathWeightingFunction}
        \tilde g_r &\coloneqq g\bigl(\trt_{\wp_r}(t)\bigr),
        &
        g_r &\coloneqq \tfrac{\tilde g_r}{\sum_{j=1}^\ell \tilde g_j}
    \end{align}
    and set the $c$\nbhy component of $\vec p_{s,\hat s}$ to $\sum_{r\in \{1,\ldots,\ell\}\colon \hat s\text{ is first street in }\wp_r} g_r$.
\end{definition}
Effectively, the weighting function~$g$ in \cref{eq:pathWeightingFunction} assigns some kind of quality to each path. These qualities are rescaled to sum up to one. Since routing does not assign shares of traffic flow to paths, but to roads, the routing~$\vec p_{s,\hat s}$ is then a sum over the qualities of all paths using road~$\hat s$ first. A typical choice for the function~$g$ could be $\trt\mapsto\trt^{-c}$ or $\trt\mapsto\exp(-c\trt)$ with $c\in\IR_{\geq0}$. However, any decreasing function might be chosen.

In the numerics section \cref{sec:examples} we mostly use $k$~shortest path routing from \cref{def:kShortestPath}. However, we also employ routing based on the future evolution of traffic (optimal control), see \cref{ssec:routingOptimization}

\begin{remark}[Problem in full complexity]
    Formulated in its full complexity, the problem reads for each street~$s$ with $t\in(0,T)$ and $x\in(0,L_s)$
    \begin{align}
        \vec q_s(0,x) &= \vec q_{0,s}(x),
        \\
        \vec b_s(0) &= \vec b_{0,s},
        \\
        0 &= \partial_t \vec q_s(t,x) + \partial_x\bigl(\vec q_s(t,x) v_s[q_s](t, x)\bigr),
        \\
        \vec q_s(t,0) v_s[q_s](t,0) &= \ul[\vec]_s(t) v_s[q_s](t,0),
        \\
        \dot{\vec b}_s(t) &= \sum_{r\in\zeta_s^-} \vec p_{r,s}(t) \vec q_r(t,L_r) v_r[q_r](t,L_r) - \ul[\vec]_s(t) v_s[q_s](t, 0),
        \\
        \ur_s(t) &= \qmax_s \max_{\{ r\in\zeta_s^+ \,\vert\, \exists c\colon [\vec q_s(t, L_s)]_c > 0 \land [\vec p_{s, r}(t)]_c > 0 \}} \tfrac{b_r(t)}{\bmax_r},
        \\
        \ul[\vec]_s(t) &= \Cases[,]{
            \Bigl(\frac{\vec b_s(t)}{\bmax_s}\Bigr)^{\frac{1}{m_s}} \qmax_s \tfrac{\vec b_s(t)}{b_s(t)}
        }{
            b_s(t)>0
        }{
            \boldsymbol0
        }
        \\
        v_s[q_s](t,x) &= V_s\Biggl(\tfrac{1}{\eta_s}\biggl(\int_x^{L_s} \gamma_s\bigl(\tfrac{y-x}{\eta_s}\bigr) q_s(t,y) + \ur_s(t)\int_{L_s}^\infty \gamma_s\bigl(\tfrac{y-x}{\eta_s}\bigr) \,\dd y\biggr)\Biggr)
\end{align}
for a given initial data~$\vec q_{0,s}$ and~$\vec b_{0,s}$ as in \cref{eq:initialDensity,eq:initialBufferLoad}, convolution kernels~$\gamma_s$ and velocities~$V_s$ as in \cref{def:nonlocalImpact}, and routing~$\vec p$ as specified in \cref{def:routing}.
\end{remark}
\section{Notes on the implementation} \label{sec:notes}
An effective way of implementing the aforementioned model is by using a global time grid and a spatial discretization of each street, assuming constant density between two neighboring grid points. By this, the nonlocal impact collapses to a weighted sum. Given the nonlocal impact, we can determine the velocity at each grid point. Followig the method of characteristics, the velocity is used to compute the position of the grid points at the next time step, while introducing a new characteristic at location $x=0$.

Given the assumption of constant density between neighboring characteristics, we can exactly determine the mass flowing from a street into the buffer. The buffer outflux translates into a mass between the characteristic at $x=0$ and the one that starts at $x=0$ in the next time step.

For travel times, it is the easiest to assume a steady state and to compute the instantaneous travel time consisting of the time needed to travel from the start to the end of a street and the buffer waiting time as
\begin{equation}
    \trt_s(t) = \tfrac{b_s(t)}{\ul_s(t) v_s[q_s](t,0)} + \int_0^{L_s} \tfrac{1}{v_s[q_s](t,x)} \,\dd x
\end{equation}
for $t\in(0,T)$, which again collapses to a sum in the discrete setting. With that, we can compute the routing as described. The $k$~shortest path algorithm in use is an adaption of~\cite{Meral2011K}, which is based on Yen's algorithm~\cite{Yen1971Finding}.

The traffic density is given by the mass between two characteristics (known from the initial density and initial characteristics distance) divided by the current distance.

The problem in its full complexity is solved by \cref{alg:networkSolver}.
\begin{algorithm}
    \caption{Network solver.}
    \label{alg:networkSolver}
    \begin{algorithmic}[1]
        \Require{network structure, \ie, streets, crossroads, commodities with their respective descriptors and connectivity information}
        \Require{time grid~$(t_i)_{i=0,\ldots,\NT}$}
        \Require{initial characteristics~$\xi_{s,j,0}$}
        \Require{initial densities~$\vec q_{s,j,0}$}
        \Require{initial buffer loads~$\vec b_{s,0}$}
        
        \ForAll{street~$s$} \Comment{Initialization.}
            \State $\vec b_s(t_0) \coloneqq \vec b_{s,0}$
            \ForAll{characteristics $j$}
                \State $\xi_{s,j}(t_0) = \xi_{s,j,0}$
                \State $\vec q_s(t_0, \xi_{s,j,0}) \coloneqq \vec q_{s,j,0}$
            \EndFor
        \EndFor
        
        \For{$i=1,\ldots,\NT$} \Comment{Time stepping.}
            \ForAll{street~$s$}
                \State compute travel time~$\trt_s(t_{i-1})$
            \EndFor
            \ForAll{street~$s$}
                \ForAll{street~$r$}
                    \State compute routing~$\vec p_{s,r}(t_{i-1})$
                \EndFor
            \EndFor
            \ForAll{street~$s$}
                \ForAll{characteristics~$j$}
                    \State $\xi_{s,j}(t_i) \coloneqq \xi_{s,j}(t_{i-1}) + (t_i-t_{i-1}) v_s[q_s](t_{i-1}, \xi_{s,j}(t_{i-1}))$
                \EndFor
                \State $\influx[\vec]_s(t_{i-1}) \coloneqq \sum_{r\in\zeta_s^+} \vec p_{r,s}(t_{i-1}) \vec q_r(t_{i-1}, L) v_r[q_r](t_{i-1}, L)$
                \State $\outflux[\vec]_s(t_{i-1}) \coloneqq \vec \ul_s(\vec b_s(t_{i-1})) v_s[q_s](t_{i-1},0)$
                \State $\vec b_s(t_i) \coloneqq \vec b_s(t_{i-1}) + (t_i-t_{i-1}) (\influx[\vec]_s(t_{i-1})-\outflux[\vec]_s(t_{i-1})$
                \ForAll{characteristics~$j$}
                    \State $\vec q_s(t_i, \xi_{s,j}(t_i)) \coloneqq \vec q_{s,j,0} \frac{\xi_{s,j+1,0}-\xi_{s,j,0}}{\xi_{s,j+1}(t_i)-\xi_{s,j}(t_i)}$
                \EndFor
            \EndFor
        \EndFor
    \end{algorithmic}
\end{algorithm}
\section{Examples of use} \label{sec:examples}
In the following, we would like to show some examples to demonstrate the capabilities of the model and its implementation. We employ \texttt{MATLAB} R2025a for the showcases.
\subsection{Preliminaries}
Before starting into the examples, we need to give some preliminary explanations on the visualization.
\subsubsection{Network structure} In network visualizations as \cref{fig:multiRouteGraph:network,fig:latticeGraph:network,fig:routingOptimization:network}, cross sections are visualized as circular nodes containing their indices. Streets may also be indexed. Arrow heads represent a street's flow direction. Networks are visualized as digraphs. Edges may be bound to either not overlap when there is another edge with opposite orientation, or to have a length proportional to the street length~$L_s$.
\subsubsection{Routing plots} Routing plots as \cref{fig:multiRouteGraph:routing,fig:latticeGraph:routing,fig:routingOptimization:routing} display the routing for chosen junctions and commodities over time. As long as routing is active, \ie, there is a positive density at the end of the street, solid lines are used. Dotted lines indicate inactive routing, \ie, routing that would be applied if there were drivers to be routed.
\subsubsection{Relative density plots} Density plots \cref{fig:latticeGraph:reldens,fig:multiRouteGraph:reldens,fig:routingOptimization:reldens} visualize the relative density~$\frac{q_s(t,x)}{\qmax_s}$ (\ie, summed over all commodities and divided by the street capacity) on each street based on the mass between and distance of two neighboring characteristics using a colormap. Pie charts in the beginning of each street visualize the buffer load~$\vec b(t)$, where each color corresponds to a commodity.
\subsubsection{Colormap} 
\begin{figure}
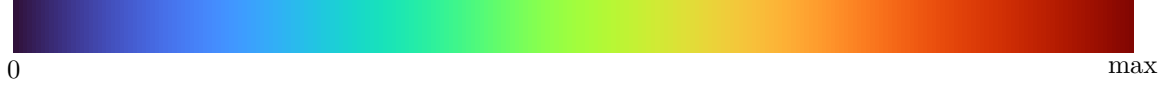

    \def\wdth{0.9\linewidth}
    \loadandscalepgfplot{1}{fig_colorbar.tex}
    \caption{Colorbar.}
    \label{fig:colorbar}
\end{figure}

All three-dimensional plots employ the same colormap as visualized in \cref{fig:colorbar}.
\subsection{Multi-route graph}
\begin{figure}
    \begin{tikzpicture}[scale=1.7]

    \foreach \x/\y [count=\i] in {-1.5/0,-1/0,-.5/-sin(60),cos(45)/sin(45),0/0,cos(15)/sin(15)}{
        \node[network node] (cr\i) at ({\x},{\y}) {\i};
    }
    \draw[network edge               ] (cr1) to node[midway, anchor=south]{\footnotesize1} (cr2);
    \draw[network edge               ] (cr2) to node[midway, anchor=east ]{\footnotesize2} (cr3);
    \draw[network edge, bend left= 28] (cr2) to node[midway, anchor=north]{\footnotesize3} (cr4);
    \draw[network edge               ] (cr2) to node[midway, anchor=south]{\footnotesize4} (cr5);
    \draw[network edge               ] (cr3) to node[midway, anchor=west ]{\footnotesize5} (cr5);
    \draw[network edge, bend left=-28] (cr3) to node[midway, anchor=west ]{\footnotesize6} (cr6);
    \draw[network edge               ] (cr5) to node[midway, anchor=south]{\footnotesize7} (cr4);
    \draw[network edge               ] (cr5) to node[midway, anchor=north]{\footnotesize8} (cr6);
\end{tikzpicture}

    \caption{Layout of the multi-route graph example. The two commodities start at node~1 and want to reach nodes~4 and~6, respectively. Edge lengths are chosen to correspond to street lengths. Arrow heads indicate the street orientation.}
    \label{fig:multiRouteGraph:network}
\end{figure}
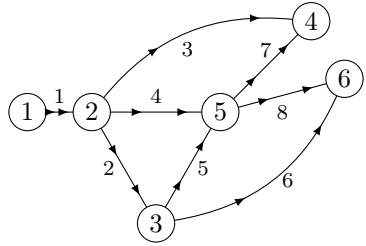%
\begin{table}
    \begin{tabular}{ccS[table-format=1.1]ccc}
        \toprule
        \thead{street $s$} & \thead{start/end node} & {\thead{length $L_s$}} & \thead{max. density $\qmax_s$} & \thead{buffer size $\bmax_s$} & \thead{initial buffer $\vec b_{s,0}$} \\ \midrule
        1 & 1/2 & 0.5 & 2 & 4 & $\begin{bmatrix} 2 & 2 \end{bmatrix}^\tp$ \\
        2 & 2/3 & 1.0 & 1 & 1 & $\begin{bmatrix} 0 & 0 \end{bmatrix}^\tp$ \\
        3 & 2/4 & 2.0 & 1 & 1 & $\begin{bmatrix} 0 & 0 \end{bmatrix}^\tp$ \\
        4 & 2/5 & 1.0 & 1 & 1 & $\begin{bmatrix} 0 & 0 \end{bmatrix}^\tp$ \\
        5 & 3/5 & 1.0 & 1 & 1 & $\begin{bmatrix} 0 & 0 \end{bmatrix}^\tp$ \\
        6 & 3/6 & 2.0 & 1 & 1 & $\begin{bmatrix} 0 & 0 \end{bmatrix}^\tp$ \\
        7 & 5/4 & 1.0 & 1 & 1 & $\begin{bmatrix} 0 & 0 \end{bmatrix}^\tp$ \\
        8 & 5/6 & 1.0 & 1 & 1 & $\begin{bmatrix} 0 & 0 \end{bmatrix}^\tp$ \\ \bottomrule
    \end{tabular}
    \caption{Selection of street parameters in the multi-route graph example.}
    \label{tab:multiRouteGraph:parameters}
\end{table}%
\begin{table}
    \begin{tabular}{cccS[table-format=1.1]}
        \toprule
        \thead{pathname} & \thead{involved nodes} & \thead{involved streets} & {\thead{path length}} \\ \midrule
        $\wp_{1,1}$      & 1, 2, 4                & 1, 3                     & 2.5                   \\
        $\wp_{1,2}$      & 1, 2, 5, 4             & 1, 4, 7                  & 2.5                   \\
        $\wp_{1,3}$      & 1, 2, 3, 5, 4          & 1, 2, 5, 7               & 3.5                   \\ \midrule
        $\wp_{2,1}$      & 1, 2, 5, 6             & 1, 4, 8                  & 2.5                   \\
        $\wp_{2,2}$      & 1, 2, 3, 5, 6          & 1, 2, 5, 8               & 3.5                   \\
        $\wp_{2,3}$      & 1, 2, 3, 6             & 1, 2, 6                  & 3.5                   \\ \bottomrule
    \end{tabular}
    \caption{Possible paths to the commodities' destinations starting in node~1. Paths~$\wp_{1,\cdot}$ belong to commodity~1, paths~$\wp_{2,\cdot}$ to commodity~2.}
    \label{tab:multiRouteGraph:paths}
\end{table}%
\begin{figure}
    \def\oid{20260612-001355-545}
    \foreach \street/\com [count=\i] in {1/1,1/2,2/2}{%
        \begin{subfigure}{0.45\linewidth}
            \loadandscalepgfplot{0.7}{fig_showRouting_ss.tex}
            \caption{Routing for commodity~\com{} on street~\street.}
        \end{subfigure}%
        \ifodd\i\relax\else\\\fi%
    }%
    \caption{Instantaneous $k$~shortest path routing for the two commodities in the multi-route graph example. When routing is inactive (\ie, there is no density to route at the end of the street), it is drawn dotted, otherwise solid.}
    \label{fig:multiRouteGraph:routing}
\end{figure}
\begin{figure}
    \def\oid{20260612-001355-545}
    \foreach\i in {0,...,8}{%
        \includegraphics[trim=50 50 50 50,clip,width=0.25\linewidth]{fig_col_\oid_t\i.000_RELDENS.png}\llap{\tiny$t=\i.0$}\allowbreak%
        \includegraphics[trim=50 50 50 50,clip,width=0.25\linewidth]{fig_col_\oid_t\i.500_RELDENS.png}\llap{\tiny$t=\i.5$}\allowbreak%
    }
    \caption{Relative densities on streets of the multi-route graph every 0.5~time units, from left to right, from top to bottom.}
    \label{fig:multiRouteGraph:reldens}
\end{figure}%
In the first example we want to develop an idea of $k$~shortest path routing and the traffic flow problem in general.

Drivers of commodity~1 with destination~4 and of commodity~2 with destination~6 both start at node~1. Their respective paths can be seen in \cref{tab:multiRouteGraph:paths} and confirmed by \cref{fig:multiRouteGraph:network}. The paths lengths are found by adding the respective street lengths displayed in \cref{tab:multiRouteGraph:parameters}. Common street parameters are velocity $V_s(w)=1-w$, look-ahead distance $\eta_s=1$, constant kernel $\gamma_s(y)=\chi_{[0,1]}(y)$, initial density $\vec q_{0,s}\equiv\vec 0$, and outflux exponent $m_s=4$. We apply $k$~shortest path routing with number of paths $k=3$ (thus taking all potential routs into consideration) and path weighting $g(\trt)=\trt^{-15}$. The scenario is solved with time step size \num{0.025} until the network evacuates.

Looking at the first commodity, paths~$\wp_{1,1}$ and~$\wp_{1,2}$ are equally long, whereas path~$\wp_{1,3}$ is longer, therefore we expect the latter to be rarely used. On first glance, we could expect paths~$\wp_{1,1}$ and~$\wp_{1,2}$ to be utilized equally, but $\wp_{1,2}$ relies on street~4, which is also used by the second commodity on path~$\wp_{2,1}$. Consequently, we expect that quickly the majority uses~$\wp_{1,1}$, a bit less~$\wp_{1,2}$, and only a few drivers consider~$\wp_{1,3}$.
The second commodity, however, wants to primarily use path~$\wp_{2,1}$ since this is the shortest one, but this one will be congested quickly as both commodities want to use it, forcing some drivers to switch to the remaining routes. As street~8 is shared by routes~$\wp_{2,1}$ and~$\wp_{2,2}$, drivers will then prefer to take~$\wp_{2,3}$ in order to have street~6 for themselves, allowing them to drive faster. Note that members of the second commodity are not routed onto street~3, because they cannot reach their destination from there.
These are exactly the effects that can be seen in \cref{fig:multiRouteGraph:routing}. \Cref{fig:multiRouteGraph:reldens} supports that by showing the relative density at different points in time.
\subsection{Lattice graph}
\begin{figure}
    \begin{subfigure}[t]{0.7\linewidth}
        \begin{tikzpicture}
    \foreach \y in {0,...,6}{
        \foreach \x in {0,...,6}{
            \node[network node] (cr\x\y) at (\x,\y) {\pgfmathparse{int(7*\y+\x+1)}\pgfmathresult};
        }
    }
    \node[network node] (cr-13) at (-1, 3) {50};
    \node[network node] (cr3-1) at ( 3,-1) {51};
    \node[network node] (cr73)  at ( 7, 3) {52};
    \node[network node] (cr37)  at ( 3, 7) {53};
    \foreach \x in {0,...,5}{
        \pgfmathsetmacro\xx{int(\x+1)}
        \foreach \y in {0,...,6}{
            \ifnum\y=3\relax
                \def\col{MatlabBlue}
            \else
                \def\col{}
            \fi
            \draw[network edge, bend left=15, \col] (cr\x\y) to (cr\xx\y);
            \draw[network edge, bend left=15, \col] (cr\xx\y) to (cr\x\y);
            \draw[network edge, bend left=15, \col] (cr\y\x) to (cr\y\xx);
            \draw[network edge, bend left=15, \col] (cr\y\xx) to (cr\y\x);
        }
    }
    \draw[network edge, MatlabRed   ] (cr-13) to (cr03);
    \draw[network edge, MatlabRed   ] (cr3-1) to (cr30);
    \draw[network edge, MatlabPurple] (cr63)  to (cr73);
    \draw[network edge, MatlabPurple] (cr36)  to (cr37);
\end{tikzpicture}
        \caption{Full network (without edge numbering).}
    \end{subfigure}%
    \begin{subfigure}[t]{0.3\linewidth}
        \begin{tikzpicture}[scale=1.4]
    \node[network node] (cr51) at ( 0,-1) {51};
    \node[network node] (cr3)  at (-1, 0) {3};
    \node[network node] (cr4)  at ( 0, 0) {4};
    \node[network node] (cr5)  at ( 1, 0) {5};
    \node[network node] (cr11) at ( 0, 1) {11};
    \draw[network edge, bend left=15            ] (cr3)  to node[midway, anchor=south] {\footnotesize  7} (cr4);
    \draw[network edge, bend left=15            ] (cr4)  to node[midway, anchor=north] {\footnotesize  9} (cr3);
    \draw[network edge, bend left=15            ] (cr4)  to node[midway, anchor=south] {\footnotesize 10} (cr5);
    \draw[network edge, bend left=15            ] (cr5)  to node[midway, anchor=north] {\footnotesize 12} (cr4);
    \draw[network edge, bend left=15, MatlabBlue] (cr4)  to node[midway, anchor=east ] {\footnotesize 11} (cr11);
    \draw[network edge, bend left=15, MatlabBlue] (cr11) to node[midway, anchor=west ] {\footnotesize 33} (cr4);
    \draw[network edge,               MatlabRed ] (cr51) to node[midway, anchor=west ] {\footnotesize170} (cr4);
\end{tikzpicture}
        \caption{Neighborhood of node~4 (with edge numbering).}
    \end{subfigure}
    \caption{Layout of the lattice graph example. The two commodities start at nodes~50 and~51 and want to reach~52 and~53, respectively. Edge lengths are chosen to correspond to street lengths. Arrow heads indicate the street orientation. Most edges are not numbered to avoid visual clutter.}
    \label{fig:latticeGraph:network}
\end{figure}
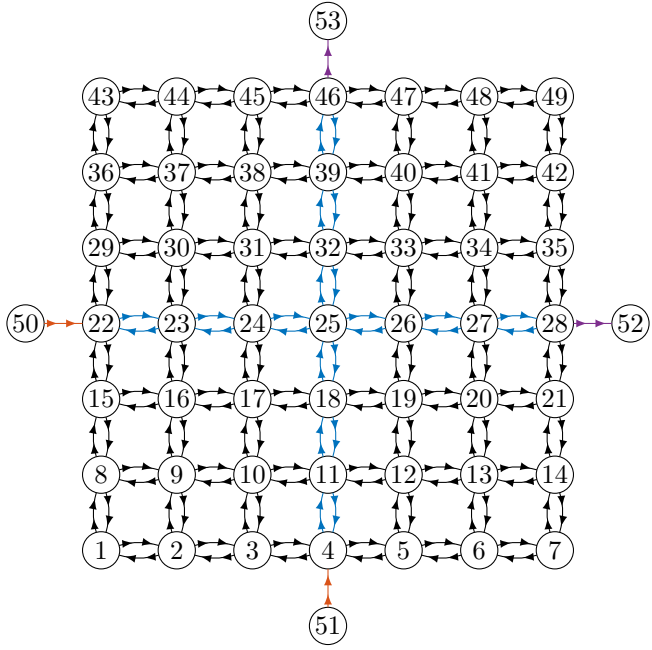
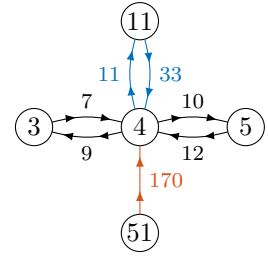
\begin{table}
    \begin{tabular}{cS[table-format=1.1]S[table-format=2.1]c}
        \toprule
        \thead{street color} & {\thead{max.\@ density~$\qmax_s$}} & {\thead{buffer size~$\bmax_s$}} & \thead{description} \\
        \midrule
        \textcolor{MatlabBlue}{blue}     & 2 &  0.5 & main street   \\
        \textcolor{MatlabRed}{red}       & 2 & 15.0 & access street \\
        \textcolor{MatlabPurple}{purple} & 2 &  0.5 & exit street   \\
        black                            & 1 &  0.5 & side street   \\
        \bottomrule
    \end{tabular}
    \caption{Selection of street parameters in the lattice graph example. For colors refer to \cref{fig:latticeGraph:network}.}
    \label{tab:latticeGraph:parameters}
\end{table}
\begin{figure}
    \def\oid{20260612-001443-517}
    \def\street{170}
    \def\com{2}
    \loadandscalepgfplot{1}{fig_showRouting_ss.tex}
    \caption{Routing for commodity~2 on street~170 (the access street) in the lattice graph example, cf.\@ \cref{fig:latticeGraph:network}. Street $r=9$ leads to the left, street $r=10$ to the right, and street $r=11$ points upwards. When routing is inactive (\ie, there is no density to route at the end of the street), it is drawn dotted, otherwise solid.}
    \label{fig:latticeGraph:routing}
\end{figure}
\begin{figure}
    \def\oid{20260612-001443-517}
    \foreach\i in {0,...,34}{%
        \includegraphics[width=0.14\linewidth]{fig_col_\oid_t\i.000_RELDENS.png}%
        \llap{\tiny$t=\i.0$}\allowbreak%
        \includegraphics[width=0.14\linewidth]{fig_col_\oid_t\i.500_RELDENS.png}%
        \llap{\tiny$t=\i.5$}\allowbreak%
    }
    \caption{Relative densities on streets of the lattice graph every 0.5~time units, from left to right, from top to bottom.}
    \label{fig:latticeGraph:reldens}
\end{figure}
The second example showcases that the code is indeed able to handle nontrivial networks of streets like the New York-like fully symmetric network in \cref{fig:latticeGraph:network} consisting of two commodities, 53~crossroads and 172~streets. Some of its parameters are collected in \cref{tab:latticeGraph:parameters}. The scenario is solved with time step size \num{0.025} until the network evacuates. We apply $k$~shortest path routing with $k=100$ and $g(\trt)=\trt^{-5}$. We further choose
\begin{align}
    V_s(w) &= 1-w,
    &
    \eta_s &= 1,
    &
    \gamma_s(y) &= \chi_{[0,1]}(y),
    &
    L_s &= 1,
    &
    m_s &= 4
\end{align}
for all streets~$s$. The initial density is negligible with $\vec q_{0,s}\equiv\begin{bmatrix}0&0\end{bmatrix}^\tp$. The buffer of the access street starting in node~50 is initially filled with drivers of commodity~1, \ie, its initial buffer load is $\begin{bmatrix}15&0\end{bmatrix}^\tp$, whereas the access street~170 starting at node~51 contains only the second commodity, $\vec b_{0,170}=\begin{bmatrix}0&15\end{bmatrix}^\tp$. All other buffers are initially empty.

A certain amount of drivers belonging to the first commodity starts at the very left, node~50, and wants to reach the very right, node~52. The other commodity, however, wants to reach the very top, node~53, starting from the bottom, node~51. Apparently, the most efficient path is to completely ignore the side streets (\eg, residential areas) and to focus on the main streets, since the latter form a direct connection from start to destination and have a higher capacity. Consequently, many drivers choose this route, leading to an increase of travel time, making the main street less attractive, causing more people to accept the detour caused by leaving the main street. This can be seen in \cref{fig:latticeGraph:routing,fig:latticeGraph:reldens}. In general, the routing at each street can be retrieved, but is not shown.

The zig-zag routing in the first seconds of the simulation is due to $k$~shortest path routing. After some time, the network is closer to an equilibrium state, meaning that the travel times on all used paths are close, thus leading to a more regular routing behavior. Around $t=20$, the buffers at the entry nodes become empty, and the network runs dry over time.

The unoptimized computation, performed on a single core of a high-end desktop computer equipped with an Intel Core i9\nbhy14900KF processor, took approximately \num{4.3}~hours for \num{2000} time steps, \ie, in total approximately \num{350000} characteristics.
\subsection{Routing optimization in time-dynamic Braess paradox} \label{ssec:routingOptimization}
The final example focuses on the minimization of the travel times under idealized conditions, \ie, we know everything, especially at any time the future behavior of the traffic state on the network for a specified routing~\cite{Bayen2019Time}. This leads to an optimal control problem~\cite{Ancona2021optimization} over the entire time horizon with an objective function yet to be specified.

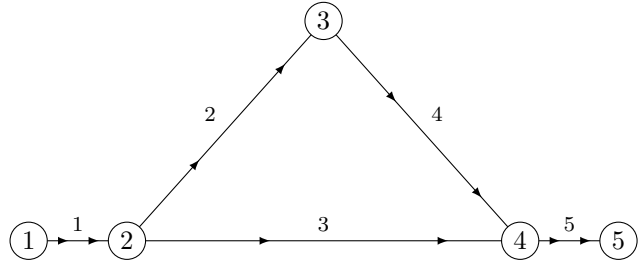
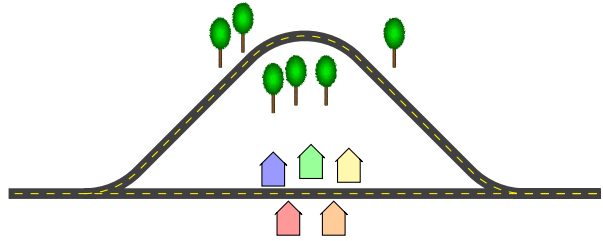
\begin{figure}
    \begin{subfigure}{\linewidth}
        \begin{tikzpicture}[scale=1.3]
    \foreach \x/\y [count=\i] in {0/0,1/0,3/sqrt(5),5/0,6/0}{
        \node[network node] (cr\i) at ({\x},{\y}) {\i};
    }
    \draw[network edge] (cr1) to node[midway, anchor=south     ] {\footnotesize1} (cr2);
    \draw[network edge] (cr2) to node[midway, anchor=south east] {\footnotesize2} (cr3);
    \draw[network edge] (cr2) to node[midway, anchor=south     ] {\footnotesize3} (cr4);
    \draw[network edge] (cr3) to node[midway, anchor=south west] {\footnotesize4} (cr4);
    \draw[network edge] (cr4) to node[midway, anchor=south     ] {\footnotesize5} (cr5);
\end{tikzpicture}
        \caption{The single commodity with destination node~5 starts in node~1. Edge lengths are chosen to correspond to street lengths. Arrow heads indicate the street orientation.}
        \label{fig:routingOptimization:network}
    \end{subfigure}\\
    \begin{subfigure}{\linewidth}
         \begin{tikzpicture}[scale=1]
    \def\mypath{(0,0) arc(270:315:1) -- coordinate[pos=0.2] (hlpr1) ++(1.5,1.5) arc(135:45:1) -- ++(1.5,-1.5) coordinate[pos=0.8] (hlpr2) arc(225:270:1) -- ++(1,0) -- (-1,0)}
    \draw[darkgray, line width=4] \mypath;
    \draw[yellow, dashed] \mypath;
    \path
        (2.5, 0.40) pic[transform shape, scale=0.3, fill=blue!40!white  ]{house}
        (3.0, 0.50) pic[transform shape, scale=0.3, fill=green!40!white ]{house}
        (3.5, 0.45) pic[transform shape, scale=0.3, fill=yellow!40!white]{house}
        (3.3,-0.25) pic[transform shape, scale=0.3, fill=orange!40!white]{house}
        (2.7,-0.25) pic[transform shape, scale=0.3, fill=red!40!white   ]{house};
    \path
        (1.8,2.0) node {\scalebox{0.15}{\tikz\pic{tree};}}
        (2.1,2.2) node {\scalebox{0.15}{\tikz\pic{tree};}}
        (2.5,1.4) node {\scalebox{0.15}{\tikz\pic{tree};}}
        (2.8,1.5) node {\scalebox{0.15}{\tikz\pic{tree};}}
        (3.2,1.5) node {\scalebox{0.15}{\tikz\pic{tree};}}
        (4.1,2.0) node {\scalebox{0.15}{\tikz\pic{tree};}};
\end{tikzpicture}
        \caption{A short street leads through the city, but has unfavorable conditions at the junction behind opposed to the slightly longer path around the city. The tree pics are adapted from~\cite{Kate2014Drawing}.}
        \label{fig:routingOptimization:sketch}
    \end{subfigure}
    \caption{Layout and sketch of the routing optimization example.}
\end{figure}
\begin{table}
    \begin{tabular}{cccS[table-format=1.1]cc}
        \toprule
        \thead{street $s$} & \thead{start/end node} & \thead{length $L_s$} & {\thead{max. density $\qmax_s$}} & \thead{buffer size $\bmax_s$} & \thead{initial buffer $\vec b_{s,0}$} \\ \midrule
        1 & 1/2 & 1 & 1.0 & 3 & 3 \\
        2 & 2/3 & 3 & 1.0 & 1 & 0 \\
        3 & 2/4 & 4 & 1.0 & 1 & 0 \\
        4 & 3/4 & 3 & 1.0 & 1 & 0 \\
        5 & 4/5 & 1 & 0.5 & 1 & 0 \\ \bottomrule
    \end{tabular}
    \caption{Selection of street parameters in the routing optimization example.}
    \label{tab:routingOptimization:parameters}
\end{table}
\begin{figure}
    \def\street{1}%
    \def\com{1}%
    \begin{subfigure}{0.5\linewidth}
        \def\oid{20260612-001510-494}%
        \loadandscalepgfplot{0.7}{fig_showRouting_ss.tex}
        \caption{$k$~shortest path routing.}
        \label{fig:routingOptimization:routing:a}
    \end{subfigure}%
    \begin{subfigure}{0.5\linewidth}
        \def\oid{20260612-004054-125}%
        \loadandscalepgfplot{0.7}{fig_showRouting_ss.tex}
        \caption{Optimal constant routing.}
        \label{fig:routingOptimization:routing:b}
    \end{subfigure}%
    \\%
    \begin{subfigure}{0.5\linewidth}
        \def\oid{20260612-115421-965}%
        \loadandscalepgfplot{0.7}{fig_showRouting_ss.tex}
        \caption{Optimal time-dependent routing.}
        \label{fig:routingOptimization:routing:c}
    \end{subfigure}%
    \begin{subfigure}{0.5\linewidth}
        \def\oid{20260612-001510-494}
        \loadandscalepgfplot{0.7}{fig_showTotalTTOverTime.tex}
        \caption{Integrand of total travel time~$\trt^\mathrm{tot}$.}
        \label{fig:routingOptimization:routing:d}
    \end{subfigure}
    \caption{\subref{fig:routingOptimization:routing:a}--\subref{fig:routingOptimization:routing:c}: Routing on street~1 in the routing optimization example. When routing is inactive (\ie, there is no density to route at the end of the street) it is drawn dotted, otherwise solid. \subref{fig:routingOptimization:routing:d}: Integrand of the total travel time~$\trt^\mathrm{tot}$.}
    \label{fig:routingOptimization:routing}
\end{figure}
\begin{figure}
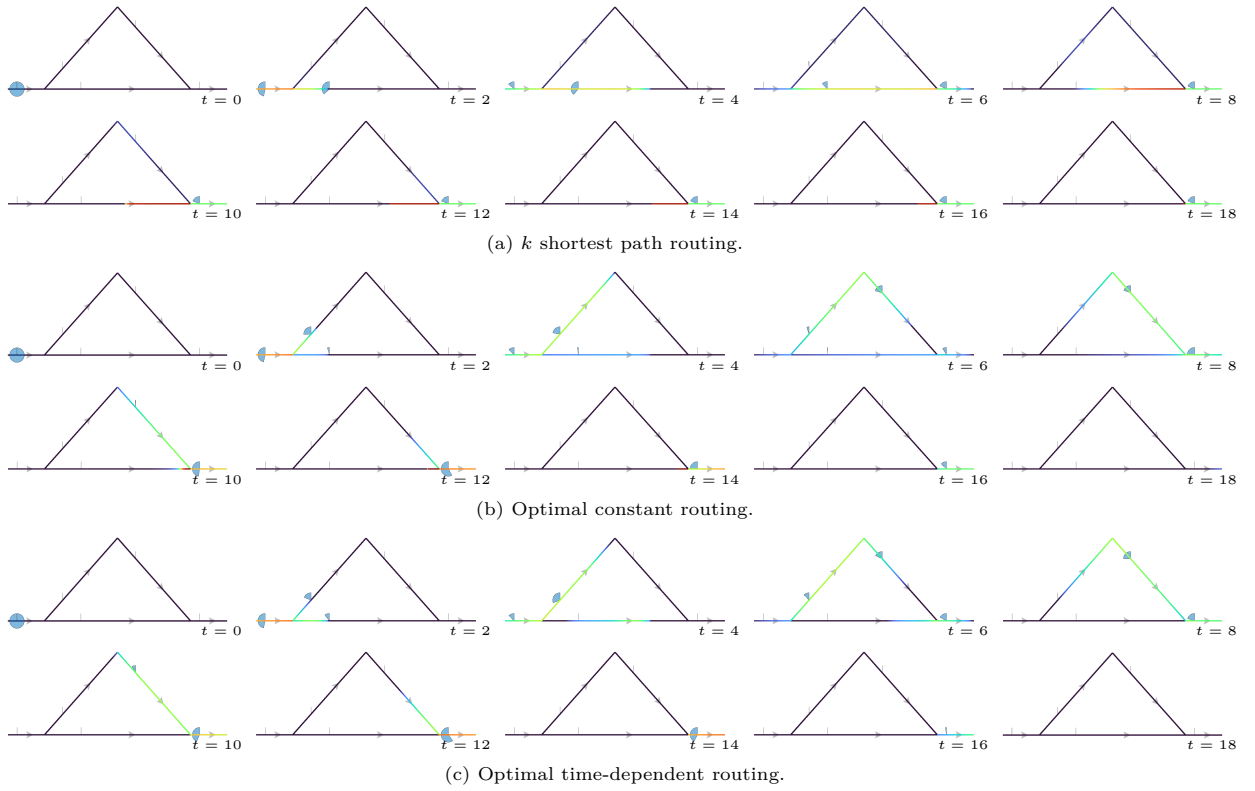

    \begin{subfigure}{\linewidth}
        \def\oid{20260612-001510-494}
        \foreach\i in {0,2,...,18}{%
            \includegraphics[width=0.2\linewidth]{fig_col_\oid_t\i.000_RELDENS.png}%
            \llap{\tiny$t=\i$}\allowbreak%
        }
        \caption{$k$~shortest path routing.}
        \label{fig:routingOptimization:reldens:a}
    \end{subfigure}\\
    \begin{subfigure}{\linewidth}
        \def\oid{20260612-004054-125}
        \foreach\i in {0,2,...,18}{%
            \includegraphics[width=0.2\linewidth]{fig_col_\oid_t\i.000_RELDENS.png}%
            \llap{\tiny$t=\i$}\allowbreak%
        }
        \caption{Optimal constant routing.}
        \label{fig:routingOptimization:reldens:b}
    \end{subfigure}\\
    \begin{subfigure}{\linewidth}
        \def\oid{20260612-115421-965}
        \foreach\i in {0,2,...,18}{%
            \includegraphics[width=0.2\linewidth]{fig_col_\oid_t\i.000_RELDENS.png}%
            \llap{\tiny$t=\i$}\allowbreak%
        }
        \caption{Optimal time-dependent routing.}
        \label{fig:routingOptimization:reldens:c}
    \end{subfigure}\\
    \caption{Relative densities on streets of the routing optimization example graph every 2~time units, from left to right, from top to bottom.}
    \label{fig:routingOptimization:reldens}
\end{figure}

We have a single commodity starting in node~1 trying to reach destination node~5 (see \cref{fig:routingOptimization:network,fig:routingOptimization:sketch}), for which they can use two routes, namely~$\wp_1$ including streets~1, 3, and~5, and the second path, $\wp_2$, using streets~1, 2, 4, and~5. From \cref{tab:routingOptimization:parameters} we derive that $\wp_1$ has length~6, \ie, it is shorter than the detour in path~$\wp_2$, which has length~8. The only routing decision is to be made in node~2, whether to take the detour or not.
Common street parameters are velocity $V_s(w)=1-w^4$, look-ahead distance $\eta_s=3$, linear kernel $\gamma_s(y)=2(y-1)\chi_{[0,1]}(y)$, initial density $q_{0,s}\equiv0$, and outflux exponent $m_s=2$. The scenario is solved with time step size~\num{0.025} until the network evacuates. Particularly the choice for the velocity~$V$ allows us to have a relatively high velocity although having high density.

To measure the time that all the drivers need, we consider the total travel time
\begin{equation}\label{eq:totalTravelTime}
    \trt^\mathrm{tot} = \int_0^T \int_0^t \outflux_1(s) - \influx_5(s) \,\dd s \,\dd t.
\end{equation}
The total travel time defined in \cref{eq:totalTravelTime} is a Wasserstein\nbhy1 metric considering the cumulative mass having entered street~2 and the mass having entered the buffer of street~5, \ie, having left the relevant part of the system. We compare three different types of routing:
\begin{enumerate}
    \item instantaneous $k$~shortest path routing (in this example, $k=2$, as we only have two routes, and path weighting $g(\trt)=\trt^{-15}$) as in \cref{def:kShortestPath},
    \item optimal constant routing, \ie, we optimize for the \emphTwo{constant} share of cars routed onto path~$p_2$, and
    \item optimal time-dependent routing, \ie, we optimize the share of cars routed onto path~$p_2$ at 31 fixed time points and interpolate linearly, while slightly penalizing steep gradients in order to achieve smoother routing by adding $\alpha\lVert p_{1,2}\rVert_{L^2((0,T))}$ with $\alpha\in\IR_{>0}$ to the total travel time.
\end{enumerate}
The third approach must always be superior to the second one, \ie, yield lower total travel times, since it is more flexible. We also expect it to be better than the first one, as long as the number of routing time points is sufficiently high, since the optimization can \enquote{learn} from previous iterations and $k$~shortest path might fall for bottlenecks. Finally, if the routing time points are the same as for the instantaneous routing, the instantaneous one becomes part of the admissible set for the optimizer, so the optimization approach must be at least as good. However, with respect to computation time, the third approach is also the most costly one. Whether routing one or two is better depends on the situation, because optimal constant routing does have the possibility to descend from previous iterations, but is also restricted to constant solutions. Note that $k$~shortest path in a sense models egoistic routing decisions, while the others aim for social optima. Therefore, we can effectively determine the Price of Anarchy~\cite{Graf2023Price}.

In this example, path~$p_1$ is the shorter and therefore preferred one. However, we make street~3 very sensitive with respect to the buffer of its successor, \ie, we replace---only on this specific road---the nonlocal impact~$W_3[q_3, \gamma_3, \eta_3]$ from \cref{eq:nonlocalImpact} with
\begin{equation}\label{eq:alteredNonlocalImpact}
    (t,x)
    \mapsto
    \frac{1}{\eta_3} \biggl(
                                   \int_x^{L_3}      \gamma_3\bigl(\tfrac{y-x}{\eta_3}\bigr) q_3(t, y) \,\dd y + 
        \max\{5\ur_3(t), \qmax_3\} \int_{L_3}^\infty \gamma_3\bigl(\tfrac{y-x}{\eta_3}\bigr)           \,\dd y
    \biggr)
\end{equation}
with $t\in(0,T)$ and $x\in(0,L_3)$. Hence, as long as the buffer of street~5 is empty, route~$\wp_1$ is---intuitively---preferable and therefore used by almost all drivers in the case of instantaneous $k$~shortest path routing due to the choice of \cref{eq:alteredNonlocalImpact}. However, sending a lot of mass onto the said route leads to large buffer loads as soon as the first drivers reach the end of street~3, causing congestion and high travel times. When this happens, almost no cars are to be routed in node~2 any more, so it is too late for the routing to react  to the increase in travel time. This is the well-known Braess paradox~\cite{Braess1968Uber,Colombo2020Microscopic} and can be observed in \cref{fig:routingOptimization:reldens:a}. On the other hand, using one of the optimization approaches, we can foresee this issue and circumvent it by sending already at early stages some traffic participants to the detour~$\wp_2$. Especially in the third routing case we can see that initially a small burst of cars is sent to the short route~$\wp_1$, since they can leave the network before the buffer of street~5 becomes too full, as the path is shorter than the detour.
The respective routings and relative densities can be seen in \cref{fig:routingOptimization:routing,fig:routingOptimization:reldens:b,fig:routingOptimization:reldens:c}.

The resulting total travel times for the three approaches are \num{27.8781}, \num{23.8495}, and \num{22.6699}, respectively, rendering the two optimized approaches to be \qty{14.4505}{\percent} and \qty{18.6818}{\percent} better than the instantaneous $k$~shortest path routing. Whereas the difference is not too high in this simple example, this is likely to change with a more complex network topology or with a higher scaling factor in \cref{eq:alteredNonlocalImpact}. In addition, we cannot guarantee that the minima are global. As discussed, the optimized time-dependent routing outperforms the other two, but is also computationally the most expensive one. It can also be noted that interestingly the integrand of the total travel time~$\trt^\mathrm{tot}$, visualized in \cref{fig:routingOptimization:routing:d} in the time-dependent optimal routing seems to always always be close to the optimum within the $k$~shortest path and constant optimal routing options, and for $t>12.5$ even lower.

For the optimization of the second and third routing approach, we employ MATLAB's \verb|fmincon| optimizer with default settings except for the maximum number of function evaluations. For the second approach, we start with constant routing sending \qty{50}{\percent} to each path, and the algorithm terminates after \num{30} function evaluations. For the third approach, we initialize routing with the one found by $k$~shortest path routing, and observe termination after \num{2909} function evaluations.
\section{Conclusions}
This paper presents a general framework for modeling nonlocal traffic flow on road networks, emphasizing the role of routing in optimizing network performance. By considering nonlocal interactions in traffic dynamics, our model will give a comprehensive understanding of how traffic behaves across large-scale networks, since anticipatory human behavior (both locally in the choice of velocity and globally in routing) is modeled. The routing strategies we investigate---ranging from the $k$~shortest path algorithm to optimization methods---demonstrate the significant impact routing decisions can have on reducing congestion and improving network efficiency.

A central component of the approach is the coupling of nonlocal edge dynamics via buffers, which provide due to their mass preserving property a consistent mechanism for mass exchange at junctions. This structure allows the model to be applied beyond simplified intersection configurations and supports general graph-based road networks.

The model's flexibility in handling diverse traffic scenarios, including multiple commodities and complex routing strategies, makes it a valuable tool for future studies. We plan to validate the model with real-world data in future work. Additionally, further extensions of the presented mode such as incorporation of traffic lights, dynamic speed limits, and shoulder lane control as well as game theory approaches could incorporate dynamic adjustments to routing based on real-time traffic data and predictions \cite{Bayen2019Time,Graf2023Prediction,Graf2023Price}, improving the model's applicability to practical traffic management.


\section*{Declarations}
\paragraph*{Availability of data and materials}
The datasets used and/or analysed during the current study are available from the corresponding author on reasonable request.
\paragraph*{Competing interests}
The authors declare that they have no competing interests.
\paragraph*{Funding}
LP and FP have been supported by the DFG -- Project-ID 416229255 -- SFB 1411. AK has been supported by the DFG -- Project-ID 547096773.

\end{document}